\documentclass[11pt]{article}
\usepackage{enumerate}
\usepackage{amssymb,a4wide,latexsym,makeidx,epsfig,fleqn}
\usepackage{amsthm}
\usepackage{amsmath}
\usepackage{enumerate}
\usepackage{graphicx}
\usepackage{subfigure}
\usepackage{float}
\usepackage{caption}
\usepackage[colorlinks, linkcolor=black, anchorcolor=black, citecolor=blue]{hyperref}%
\newtheorem{theorem}{Theorem}[section]

\newtheorem{definition}[theorem]{Definition}
\newtheorem{lemma}[theorem]{Lemma}

\begin{document}
\textwidth 150mm \textheight 225mm
\title{Spectral radius of graphs of given size with forbidden subgraphs\thanks{Supported by the National Natural Science Foundation of China (No. 12271439).}}
\author{{Yuxiang Liu$^{a,b}$, Ligong Wang$^{a,b,}$\footnote{Corresponding author.}}\\
{\small $^{a}$  School of Mathematics and Statistics, Northwestern Polytechnical University,}\\{\small  Xi'an, Shanxi 710129, P.R. China.}\\ {\small $^{b}$ Xi'an-Budapest Joint Research Center for Combinatorics, Northwestern Polytechnical University,}\\{\small  Xi'an, Shanxi 710129, P.R. China.}\\
{\small E-mail: yxliumath@163.com, lgwangmath@163.com}}
\date{}
\maketitle
\begin{center}
\begin{minipage}{135mm}
\vskip 0.3cm
\begin{center}
{\small {\bf Abstract}}
\end{center}
{\small  Let $\rho(G)$ be the spectral radius of a graph $G$ with $m$ edges. Let $S_{m-k+1}^{k}$ be the graph obtained from $K_{1,m-k}$ by adding $k$ disjoint edges within its independent set. Nosal's theorem states that if $\rho(G)>\sqrt{m}$, then $G$ contains a triangle. Zhai and Shu showed that any non-bipartite graph $G$ with $m\geq26$ and   $\rho(G)\geq\rho(S_{m}^{1})>\sqrt{m-1}$ contains a quadrilateral unless $G\cong S_{m}^{1}$ [M.Q. Zhai, J.L. Shu, Discrete Math. 345 (2022) 112630]. Wang proved that if $\rho(G)\geq\sqrt{m-1}$ for a graph $G$ with size $m\geq27$, then $G$ contains a quadrilateral unless $G$ is one of four exceptional graphs [Z.W. Wang, Discrete Math. 345 (2022) 112973]. In this paper, we show that any non-bipartite graph $G$ with size $m\geq51$ and $\rho(G)\geq\rho(S_{m-1}^{2})>\sqrt{m-2}$ contains a quadrilateral unless $G$ is one of three exceptional graphs. Moreover, we show that if $\rho(G)\geq\rho(S_{\frac{m+4}{2},2}^{-})$ for a graph $G$ with even size $m\geq74$, then $G$ contains a $C_{5}^{+}$ unless $G\cong S_{\frac{m+4}{2},2}^{-}$, where $C_{t}^{+}$ denotes the graph obtained from $C_{t}$ and $C_{3}$ by identifying an edge, $S_{n,k}$ denotes the graph obtained by joining each vertex of $K_{k}$ to $n-k$ isolated vertices and $S_{n,k}^{-}$ denotes the graph obtained by deleting an edge incident to a vertex of degree two, respectively.
}
\end{minipage}
\end{center}
\section{Introduction}
Throughout this paper, all graphs considered are always undirected and simple. Let $G$ be a graph of order $n$ with vertex set $V(G)=\{v_{1}, v_{2},\ldots, v_{n}\}$ and size $m$ with edge set $E(G)=\{e_{1}, e_{2}, \ldots, e_{m}\}$. The neighborhood
of a vertex $u\in V(G)$ is denoted by $N_{G}(u)$. Let $N_{G}[u]=N_{G}(u)\cup\{u\}$, which is called the closed neighborhood of $u$. Let $d_{G}(u)$ be the degree of a vertex $u$. For the sake of simplicity, we omit all the subscripts if $G$ is clear from the context. The adjacency matrix of $G$ is an $n\times n$ matrix $A(G)$ whose $(i,j)$-entry is $1$ if $v_{i}$ is adjacent to $v_{j}$ and $0$ otherwise. The spectral radius $\rho(G)$ of $G$ is the
largest eigenvalue of its adjacency matrix $A(G)$.

Let $P_{n}, C_{n}, K_{1,n}$ and $K_{a,b}$ be the path of order $n$, the cycle of order $n$, the star graph of order $n+1$ and the complete bipartite graph with two parts of sizes $a, b$, respectively. Let $S_{n}^{k}$ be the graph obtained from $K_{1,n-1}$ by adding $k$ disjoint edges within its independent sets. Let $S_{n,k}$ be the graph obtained by joining each vertex of $K_{k}$ to $n-k$ isolated vertices. Let $S_{n,k}^{-}$ be the graph obtained from $S_{n,k}$ by deleting an edge incident to a vertex of degree two. Let $C_{t}^{+}$ be the graph obtained from $C_{t}$ and $C_{3}$ by identifying an edge.

Given a graph $F$, a graph $G$ is $F$-free if it does not contain $F$ as a subgraph. Let $\mathcal{G}(m, F)$ denote the family of $F$-free graphs with $m$ edges and without isolated vertices. A classic problem in extremal graph theory, known as Tur\'{a}n's problem, is that what the maximum number of edges in an $F$-free graph of order $n$ is. Nikiforov \cite{Ni2} posed a spectral version of Tur\'{a}n's problem as follows: what is the maximum spectral radius of an $F$-free graph of order $n$? This spectral Tur\'{a}n-type problem of graphs have received much attention in the past decades. For example, some new results were found in \cite{ChLZ,CiDT,CiFTZ,LiuBW, ZhL}. For more results on spectral extremal graph theory, we suggest the reader to see surveys \cite {ChZ,FSi,LiLF,Ni3}, and references therein. In contrast, the spectral Tur\'{a}n-type problem of graphs with given size is that what the maximum spectral radius of an $F$-free with $m$ edges is. Equivalently, what is a lower bound of $\rho(G)$ for a graph $G$ of size $m$ containing a subgraph $F$? Earliest, Nosal \cite{Nos} showed that if $\rho(G)>\sqrt{m}$ then $G$ contains a triangle, which is known well as a spectral Mantel's theorem. Very recently, Lin, Ning and Wu \cite{LNW} showed that if $\rho(G)\geq\sqrt{m-1}$ for a non-bipartite graph $G$ of size $m$, then $G$ contains a triangle unless $G\cong C_{5}$. Zhai and Shu \cite{ZhS} showed that if $\rho(G)\geq\rho(SK_{2,\frac{m-1}{2}})$ for a non-bipartite graph $G$ of size $m$, then $G$ contains a triangle unless $G\cong SK_{2,\frac{m-1}{2}}$, where $SK_{2,\frac{m-1}{2}}$ is the graph obtained from $K_{2, \frac{m-1}{2}}$ by subdividing an edge. Wang \cite{Wan} showed that if $\rho(G)\geq \sqrt{m-2}$ for a non-bipartite graph $G$ of size $m\geq26$, then $G$ contains a triangle unless $G$ is one of some exceptional graphs. For more details, one may refer to \cite{GLZH,LiPe,LinG} and references therein.

\noindent\begin{theorem}\label{th:ch-1.1.}{\rm(}$\cite{Wan}${\rm)}
Let $G$ be a non-bipartite and connected graph of size $m\geq26$. If $\rho(G)\geq\rho(S_{m}^{1})>\sqrt{m-1}$, then $G$ contains a quadrilateral unless $G\cong S_{m}^{1}$.
\end{theorem}

\noindent\begin{theorem}\label{th:ch-1.2.}{\rm(}$\cite{Wan}${\rm)}
Let $G$ be a graph of size $m\geq27$. If $\rho(G)\geq\sqrt{m-1}$, then $G$ contains a quadrilateral unless $G$ is one of these graphs (with possibly isolated vertices): $K_{1,m}, S_{m}^{1}, S_{m}^{e}$, or $K_{1,m-1}\cup P_{2}$, where $S_{m}^{e}$ is the graph obtained by attaching a pendent vertex to a pendent vertex of $K_{1,m-1}$.
\end{theorem}

\noindent\begin{theorem}\label{th:ch-1.3.} Let $G$ be a non-bipartite graph of size $m\geq51$. If $\rho(G)\geq\rho(S_{m-1}^{2})>\sqrt{m-2}$, then $G$ contains a quadrilateral unless $G$ is one of the following:
$S_{m}^{1}$, $C_{5}\bullet K_{1,m-5}$ and $S_{m-1}^{2}$, where $C_{5}\bullet K_{1,m-5}$ is the graph obtained by attaching a vertex of $C_{5}$ to the center vertex of $K_{1,m-5}$.
\end{theorem}

Recently, Li, Shu and Wei \cite{LiSW} characterized the extremal graph of odd size $m$ having the largest spectral radius in $\mathcal{G}(m, C_{4}^{+})$ and $\mathcal{G}(m, C_{5}^{+})$, respectively. We list them as follows.

\noindent\begin{theorem}\label{th:ch-1.4.}{\rm(}$\cite{LiSW}${\rm)} {\rm(}$i${\rm)} If $G\in\mathcal{G}(m,C_{4}^{+})$ and $m$($\geq8$) is odd, then $\rho(G)\leq\frac{1+\sqrt{4m-3}}{2}$ and equality holds if and only if $G\cong S_{\frac{m+3}{2},2}$;

{\rm(}$ii${\rm)} If $G\in\mathcal{G}(m,C_{5}^{+})$ and $m$($\geq22$) is odd, then $\rho(G)\leq\frac{1+\sqrt{4m-3}}{2}$ and equality holds if and only if $G\cong S_{\frac{m+3}{2},2}$.
\end{theorem}

Recently, Fang and You \cite{FaY} characterized the extremal graph of even size $m$ having the largest spectral radius in $\mathcal{G}(m, C_{4}^{+})$ in Theorem \ref{th:ch-1.5.}.

\noindent\begin{theorem}\label{th:ch-1.5.}{\rm(}$\cite{FaY}${\rm)} If $G\in\mathcal{G}(m,C_{4}^{+})$ and $m$($\geq 22$) is even, then $\rho(G)\leq\rho(S_{\frac{m+4}{2},2}^{-})$, and equality holds if and only if $G\cong S_{\frac{m+4}{2},2}^{-}$.
\end{theorem}
Motivated by Theorems \ref{th:ch-1.4.} and \ref{th:ch-1.5.}, we will characterize the extremal graph of even size $m$ having the maximum spectral radius in $\mathcal{G}(m,C_{5}^{+})$ as follows.

\noindent\begin{theorem}\label{th:ch-1.6.}  If $G\in\mathcal{G}(m, C_{5}^{+})$ and $m$($\geq74$) is even, then $\rho(G)\leq \rho(S_{\frac{m+4}{2},2}^{-})$, and equality holds if and only if $G\cong S_{\frac{m+4}{2},2}^{-}$.
\end{theorem}

\section{Preliminary}
In this section, we introduce some lemmas and notations. Let $X$ be the Perron vector of $G$ with coordinate $x_{v}$ corresponding to the vertex $v\in V(G)$ and $u^{\ast}$ be a vertex if $x_{u^{\ast}}=\max\{x_{v}|v\in V (G)\}$. Let $N_{i}(u)=\{v|v\in N(u)$, $d_{N(u)}(v)=i\}$, $N_{i}^{2}(u)=\{w|w\in N^{2}(u)$, $d_{N_{i}(u)}(w)\geq1\}$. Let $N[u]=N(u)\cup\{u\}$, $W=V(G)\setminus N[u]$. For a subset $S\subseteq V(G)$ and a vertex $v\in V(G)$, let $N_{S}(v)=N(v)\cap S$ and $d_{S}(v)=|N_{S}(v)|$. Let $G[S]$ be the subgraph of $G$ induced by $S$. Write $\rho=\rho(G)$. For two vertex subsets $S$ and $T$ of $V(G)$ (where $S\cap T$ may not be empty), let $e(T,S)$ denote the number of edges with one endpoint in $S$ and the other in $T$. $e(S,S)$ is simplified by $e(S)$.

\noindent\begin{lemma}\label{le:ch-2.1.} {\rm(}$\cite{ZhWF1}${\rm)} Let $u, v$ be two distinct vertices of a connected graph $G$, $\{v_{i}| i=1,2, \ldots, s\}\subseteq N(v)\setminus N(u)$, and $X=(x_{1}, x_{2}, \ldots, x_{n})^{T}$ be the Perron vector of $G$. Let $G^{\prime}=G-\sum_{i=1}^{s}v_{i}v+\sum_{i=1}^{s}v_{i}u$. If $x_{u}\geq x_{v}$, then $\rho(G)< \rho(G^{\prime})$.
\end{lemma}

\noindent\begin{lemma}\label{le:ch-2.2.}{\rm(}$\cite{Ni4}${\rm)} $\rho(S^{k}_{m-k+1})$ is the largest root of the polynomial $f(x)=x^{3}-x^{2}-(m-k)x+m-3k$, then $\sqrt{m-k}<\rho(S^{k}_{m-k+1})\leq\sqrt{m-k+1}$ for $1\leq k\leq \frac{m}{3}, m\geq4k^{2}+5k$.
\end{lemma}
\noindent {\bf Proof.} Since $f^{\prime}(x)>0$ for $x\geq\sqrt{m-k}$ and $f(\sqrt{m-k})=-2k<0, f(\sqrt{m-k+1})=\sqrt{m-k+1}-2k-1\geq0$ for $m\geq4k^{2}+5k$. Thus, we have $\sqrt{m-k}<\rho(S^{k}_{m-k+1})\leq\sqrt{m-k+1}$ for $1\leq k\leq \frac{m}{3}, m\geq4k^{2}+5k$, as desired. $\qedsymbol$

\noindent\begin{definition}\label{de:ch-2.3.}{\rm(}$\cite{CvRS}${\rm)} Given a graph $G$, the vertex partition $\Pi$: $V(G)=V_{1}\cup V_{2} \cup \ldots \cup V_{k}$ is said to be an equitable partition if, for each $u\in V_{i}$, $|V_{j}\cap N(u)|=b_{ij}$ is a constant depending only on $i,j$ ($1\leq i,j\leq k$). The matrix $B_{\Pi}=(b_{ij})$ is called the quotient matrix of $G$ with respect to $\Pi$.
\end{definition}

\noindent\begin{lemma}\label{le:ch-2.4.}{\rm(}$\cite{CvRS}${\rm)} Let $\Pi$: $V(G)=V_{1}\cup V_{2} \ldots \cup V_{k}$ be an equitable partition of $G$ with quotient matrix $B_{\Pi}$. Then $det(xI-B_{\Pi}) \mid det(xI-A(G))$. Furthermore, the largest eigenvalue of $B_{\Pi}$ is just the spectral radius of $G$.
\end{lemma}

Throughout this paper, the following equalities are used.\\

Since $A(G)X=\rho X$, we have

\begin{equation}\label{eq:ch-1}
\rho x_{u}=\sum_{v\in N_{0}(u)}x_{v}+\sum_{v\in N(u)\backslash N_{0}(u)}x_v.
\end{equation}
Since $\rho^{2}$ is the spectral radius of $A^{2}(G)$, we have
\begin{equation}\label{eq:ch-2}
\rho^{2}x_{u}=d(u)x_{u}+\sum_{v\in N(u)\backslash N_{0}(u)}d_{N(u)}(v)x_{v}+\sum_{w\in N^{2}(u)}d_{N(u)}(w)x_w.
\end{equation}
Combining with \eqref{eq:ch-1} and \eqref{eq:ch-2}, we have
\begin{equation}\label{eq:ch-3}
(\rho^{2}-\rho)x_{u}=d(u)x_{u}+\sum_{v\in N(u)\backslash N_{0}(u)}(d_{N(u)}(v)-1)x_{v}+\sum_{w\in N^{2}(u)}d_{N(u)}(w)x_w-\sum_{v\in N_{0}(u)}x_{v}.
\end{equation}

\section{Proof of Theorem 1.3.}
Let $G$ be a non-bipartite graph and $\rho(G)\geq\rho(S_{m-1}^{2})>\sqrt{m-2}\geq7$ for $m\geq51$. Recall that $W=V(G)\backslash N[u^{\ast}]$. Assume that $G$ contains no $C_{4}$, we have $N(u^{\ast})=N_{1}(u^{\ast})\cup N_{0}(u^{\ast})$, $N_{W}(u)\cap N_{W}(v)=\emptyset$ for any two vertices $u,v\in V(G)$, and $d_{N(u^{\ast})}(w)=1$ for any vertex $w\in N^{2}(u^{\ast})$. Let $N_{1}(u^{\ast})=\{u_{2i-1}u_{2i}| i\in 1,2,\ldots, 2e(N_{1}(u^{\ast}))\}$.

\begin{equation}\label{eq:ch-4}
\begin{split}
\rho^{2}x_{u^{\ast}}&=d(u^{\ast})x_{u^{\ast}}+\sum_{v\in N_{1}(u^{\ast})}x_{v}+\sum_{w\in N^{2}(u^{\ast})}d_{N(u^{\ast})}(w)x_w\\
&\leq d(u^{\ast})x_{u^{\ast}}+ \sum_{u_{2i-1}u_{2i}\in E(G[N_{1}(u^{\ast})])}(x_{u_{2i-1}}+x_{u_{2i}})+e(N(u^{\ast}),N^{2}(u^{\ast}))x_{u^{\ast}}.
\end{split}
\end{equation}

Since $G$ is $C_{4}$-free, we obtain any two vertices in $N(u^{\ast})$ have no common neighbors in $N^{2}(u^{\ast})$. Hence,

\begin{equation}\label{eq:ch-5}
e(W)=\frac{1}{2}\sum_{w\in W}d_{W}(w)
\geq\frac{1}{2}\sum_{w\in N^{2}(u^{\ast})}d_{W}(w)
\geq\frac{1}{2}|N^{2}(u^{\ast})|
\geq\frac{1}{2}\sum_{u\in N_{1}(u^{\ast})}d_{W}(u)
\end{equation}

For each $u_{2i-1}u_{2i}\in E(G[N_{1}(u^{\ast})])$, we have
$\rho x_{u_{2i-1}}=x_{u_{2i}}+x_{u^{\ast}}+\sum_{w\in N_{W}(u_{2i-1})}x_{w}$ and $\rho x_{u_{2i}}=x_{u_{2i-1}}+x_{u^{\ast}}+\sum_{w\in N_{W}(u_{2i})}x_{w}$. It follows that
\begin{equation}\label{eq:ch-6}
(\rho-1)(x_{u_{2i-1}}+x_{u_{2i}})\leq(2+d_{W}(u_{2i-1})+d_{W}(u_{2i}))x_{u^{\ast}}.
\end{equation}

Recall that $\rho(G)\geq\rho(S_{m-1}^{2})>\sqrt{m-2}\geq7$ for $m\geq51$. Combining with \eqref{eq:ch-5} and \eqref{eq:ch-6}, we obtain that

\begin{equation}\label{eq:ch-7}
\begin{aligned}
\sum_{u_{2i-1}u_{2i}\in E(G[N_{1}(u^{\ast})])}(x_{u_{2i-1}}+x_{u_{2i}})
&\leq\frac{1}{\rho-1}\sum_{u_{2i-1}u_{2i}\in E(G[N_{1}(u^{\ast})])}(2+d_{W}(u_{2i-1})+d_{W}(u_{2i}))x_{u^{\ast}}\\
&\leq \frac{e(N_{1}(u^{\ast}))}{3}x_{u^{\ast}}+\frac{1}{6}\sum_{u_{2i-1}u_{2i}\in E(G[N_{1}(u^{\ast})])}(d_{W}(u_{2i-1})+d_{W}(u_{2i}))x_{u^{\ast}}\\
&=\frac{e(N_{1}(u^{\ast}))}{3}x_{u^{\ast}}+\frac{1}{6}\sum_{u\in V(N_{1}(u^{\ast}))}d_{W}(u)x_{u^{\ast}}\\
&\leq\frac{e(N_{1}(u^{\ast}))+e(W)}{3}x_{u^{\ast}}.
\end{aligned}
\end{equation}

Combining with \eqref{eq:ch-4} and \eqref{eq:ch-7}, we get

\begin{equation}\label{eq:ch-8}
\begin{split}
\rho^{2}x_{u^{\ast}}&\leq(d(u^{\ast})+\frac{1}{3}(e(N_{1}(u^{\ast}))+e(W))+e(N(u^{\ast}),N^{2}(u^{\ast})))x_{u^{\ast}}\\
&=(m-\frac{2}{3}(e(N_{1}(u^{\ast}))+e(W)))x_{u^{\ast}}.
\end{split}
\end{equation}

Note that $\rho(G)\geq\rho(S_{m-1}^{2})>\sqrt{m-2}\geq7$. We get $e(N_{1}(u^{\ast}))+e(W)<3$, i.e., $e(N_{1}(u^{\ast}))+e(W)\leq2$. Since $G$ is a non-bipartite graph, we have $e(N_{1}(u^{\ast}))+e(W)\neq0$. Hence $1\leq e(N_{1}(u^{\ast}))+e(W)\leq2$. Now we consider the following two cases.

{\bf Case 1.} $e(W)+e(N_{1}(u^{\ast}))=2$.

In this case, we discuss the following three subcases.

{\bf Subcase 1.1.} $e(N_{1}(u^{\ast}))=2$.

In this case, we have $e(W)=0$. Suppose that $W\neq\emptyset$, without loss of generality, there exists a vertex $w\in W$. Since $G$ does not contain $C_{4}$, we have $d(w)=1$. Let $u\in N_{N(u^{\ast})}(w)$. Let $S_{m-1}^{2}=G-uw+u^{\ast}w$. By Lemma \ref{le:ch-2.1.}, we have $\rho(S_{m-1}^{2})>\rho(G)$, a contradiction. Thus  $W=\emptyset$ and $G^{\ast}\cong S_{m-1}^{2}$.

{\bf Subcase 1.2.} $e(N_{1}(u^{\ast}))=1$.

 In this case, we have $e(W)=1$. Let $w_{1}w_{2}\in E(G[W])$ be the unique edge. Assume that $u_{1}\in N_{N(u^{\ast})}(w_{1})\cap N_{N(u^{\ast})}(w_{2})$, then $G\cong G_{0}$ or $G\cong G_{1}$ (see Fig. 1). Note that $S_{m-1}^{2}=G_{i}-\{u_{1}w_{1}, u_{1}w_{2}\}+\{u^{\ast}w_{1}, u^{\ast}w_{2}\}$ for each $i\in \{0,1\}$. By Lemma \ref{le:ch-2.1.}, we have $\rho(S_{m-1}^{2})>\rho(G_{i})$ for each $i\in \{0,1\}$, a contradiction. Thus $N_{N(u^{\ast})}(w_{1})\cap N_{N(u^{\ast})}(w_{2})=\emptyset$. Without loss of generality, let $u_{1}\in N_{N(u^{\ast})}(w_{1})$ and $u_{2}\in N_{N(u^{\ast})}(w_{2})$. Then $G\cong G_{2}$ or $G\cong G_{3}$ (see Fig. 1). Note that $S_{m-1}^{2}=G_{i}-\{u_{1}w_{1}, u_{2}w_{2}\}+\{u^{\ast}w_{1}, u^{\ast}w_{2}\}$ for each $i\in \{2,3\}$. By Lemma \ref{le:ch-2.1.}, we have $\rho(S_{m-1}^{2})>\rho(G_{i})$ for each $i\in \{2,3\}$, a contradiction.

{\bf Subcase 1.3.} $e(W)=2$.

In this case, we obtain that $e(N_{1}(u^{\ast}))=0$ and $G$ possibly contains the following subgraphs (see Fig. 2). If $G$ contains $C_{6}$ as a subgraph, then $G$ is a bipartite graph, a contradiction. Assume that $G$ contains $C_{5}^{+}$ as a subgraph. Note that $S_{m-1}^{2}=C_{5}^{+}-\{w_{2}w_{3}, u_{1}w_{1}, u_{1}w_{2}\}+\{w_{1}u^{\ast}, w_{2}u^{\ast},w_{3}u^{\ast}\}$. By Lemma \ref{le:ch-2.1.}, we have $\rho(S_{m-1}^{2})>\rho(C_{5}^{+})$, a contradiction. Assume that $G$ contains $G_{4}$ as a subgraph. Note that $S_{m-1}^{2}=G_{4}-\{w_{1}w_{2},w_{2}w_{3}, u_{3}w_{2}\}+\{w_{1}u^{\ast}, w_{2}u^{\ast},w_{3}u^{\ast}\}$. By Lemma \ref{le:ch-2.1.}, we have $\rho(S_{m-1}^{2})>\rho(G_{4})$, a contradiction. For the rest graphs $G_{i}$ for $i\in\{5,6,7,8,9\}$, we have the similar operation and conclusion.

{\bf Case 2.}  $e(W)+e(N_{1}(u^{\ast}))=1$.

In this case, we discuss the following two subcases.

{\bf Subcase 2.1.} $e(N_{1}(u^{\ast}))=1$.

In this case, we have $e(W)=0$. Suppose that $W\neq\emptyset$, without loss of generality, there exists a vertex $w_{1}\in W$. Since $G$ does not contain $C_{4}$, we have $d(w_{1})=1$. Let $u_{1}\in N_{N_{1}(u^{\ast})}(w_{1})$. Then $G\cong G_{10}$ (see Fig. 3). By Lemma \ref{le:ch-2.4.}, $\rho(G_{10})$ is the largest roots of the equation $g(x)=0$, where $$g(x)=x^{4}-mx^{2}-2x+2m-7.$$ Since $g(\sqrt{m-2})=-2\sqrt{m-2}-3<0$ and $g^{\prime}(x)>0$ for $x\geq\sqrt{m-2}$. Thus $\sqrt{m-2}<\rho(G_{10})$. By Lemma \ref{le:ch-2.2.}, $\rho(S_{m-1}^{2})$ is the largest root of the equation $f(x)=0$, where $$f(x)=x^{3}-x^{2}-(m-2)x+m-6.$$ Let $$h(x)=g(x)-xf(x)=x^{3}-2x^{2}-(m-8)x+2m-7.$$ By calculation, $h^{\prime}(x)>0$ for $x\geq\sqrt{m-2}$ and $h(\sqrt{m-2})=6\sqrt{m-2}-3>0$ for $m\geq51$. Thus $\rho(G_{10})<\rho(S_{m-1}^{2})$, a contradiction. Thus $W=\emptyset$ and $G\cong S_{m}^{1}$. By Lemma \ref{le:ch-2.1.}, $\rho(S_{m}^{1})>\rho(S_{m-1}^{2})$, as desired.

{\bf Subcase 2.2.} $e(W)=1$.

In this case, we have $e(N_{1}(u^{\ast}))=0$ and $G\cong C_{5}\bullet K_{1,m-5}$ or $G\cong G_{11}$ or $G\cong G_{12}$ (see Fig. 3). By Lemma \ref{le:ch-2.4.}, $\rho(C_{5}\bullet K_{1,m-5})$, $\rho(G_{11})$ and $\rho(G_{12})$ are the largest roots of these equations $h_{1}(x)=0$, $h_{2}(x)=0$ and $h_{3}(x)=0$ respectively, where
\begin{equation}\label{eq:ch-9}
\begin{split}
h_{1}(x)&= x^{4}-x^{3}-(m-2)x^{2}-(m-3)x+m-5,\\
h_{2}(x)&=x^{5}+x^{4}-(m-1)x^{3}+x^{2}+(3m-15)x+3m-17,\\
h_{3}(x)&=x^{4}-x^{3}-(m-1)x^{2}-(m-4)x+2m-8.
\end{split}
\end{equation}
By Lemma \ref{le:ch-2.2.}, $\rho(S^{2}_{m-1})$ is the largest root of the equation $f(x)=0$. Thus $$h_{1}(x)-xf(x)=-(2m-9)x+m-5<0$$ and $\rho(C_{5}\bullet K_{1,m-5})> \rho(S^{2}_{m-1})$, as desired. Since $h_{2}(\sqrt{m-2})>0$ and $h_{2}^{\prime}(x)>0$ for $x>\sqrt{m-2}$. Thus $\rho(G_{11})< \rho(S^{2}_{m-1})$, a contradiction. Since $h_{3}(\sqrt{m-2})=m-6-2\sqrt{m-2}>0$ and $h_{3}^{\prime}(x)>0$ for $x>\sqrt{m-2}$. Hence, $\rho(G_{12})<\sqrt{m-2}$, a contradiction.

This completes the proof. $\blacksquare$

\begin{figure}[H]
\begin{centering}
\includegraphics[scale=0.25]{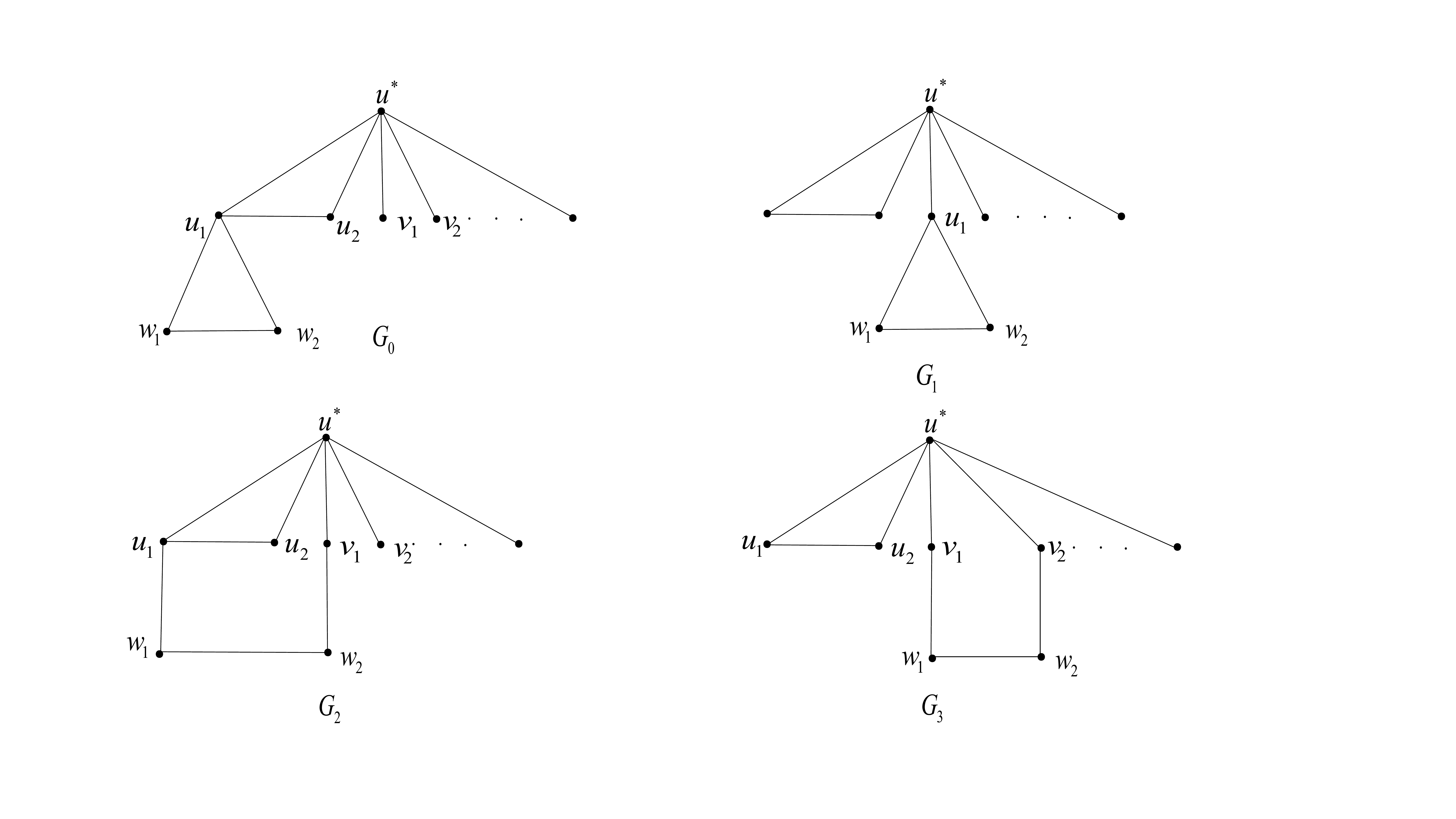}
\caption{Graphs $G_{0}-G_{3}$ of Subcase 1.2.}\label{Fig. 1.}
\end{centering}
\end{figure}

\begin{figure}[H]
\begin{centering}
 \subfigure[]{
  \includegraphics[scale=0.25]{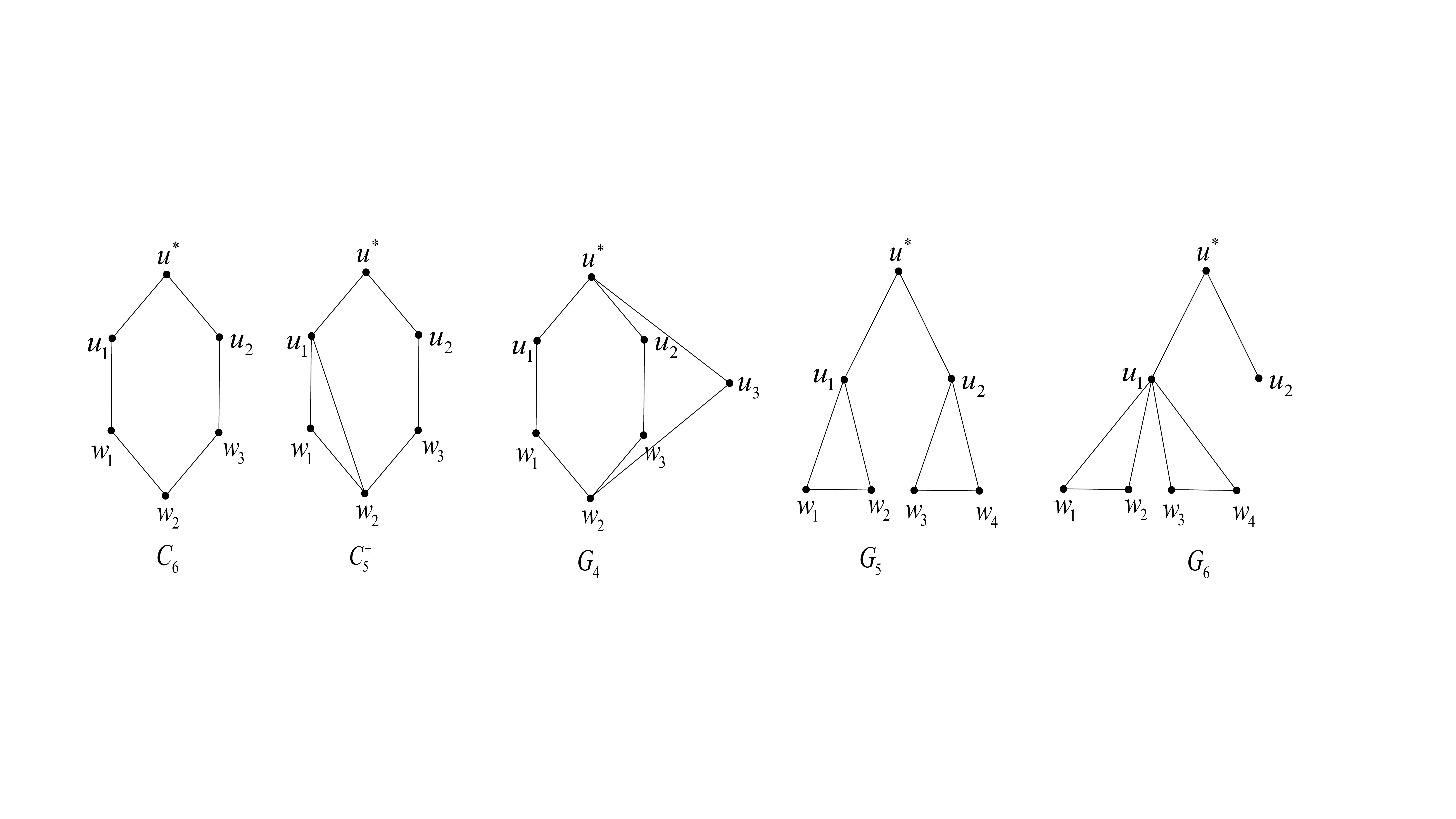}
   }
   \quad
   \subfigure[]{
       \includegraphics[scale=0.25]{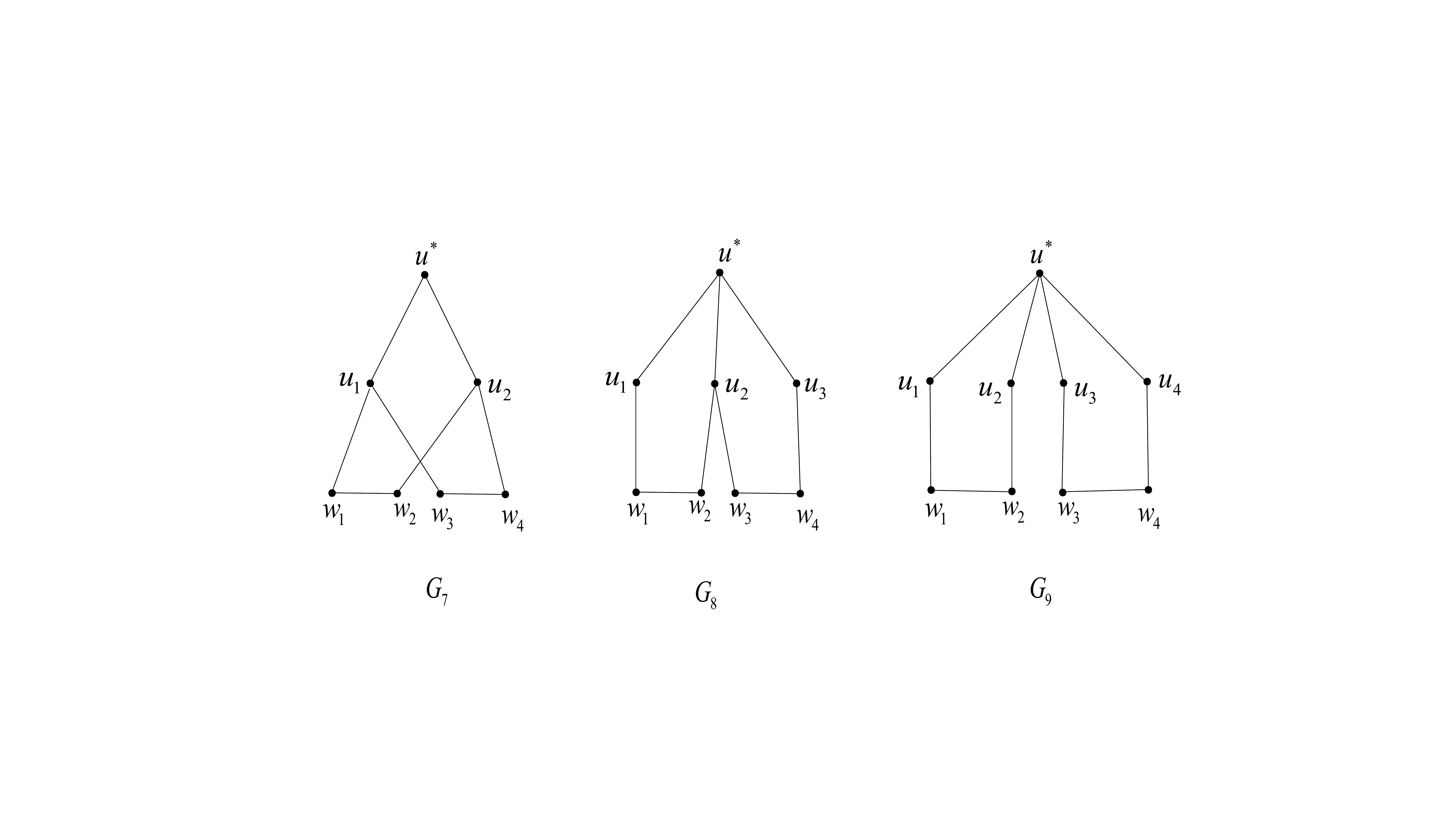}
       }
       \caption{Graphs $C_{6}, C_{5}^{+}$ and $G_{4}-G_{9}$ of Subcase 1.3.}\label{2.}
\end{centering}
\end{figure}

\begin{figure}[H]
\begin{centering}
\includegraphics[scale=0.25]{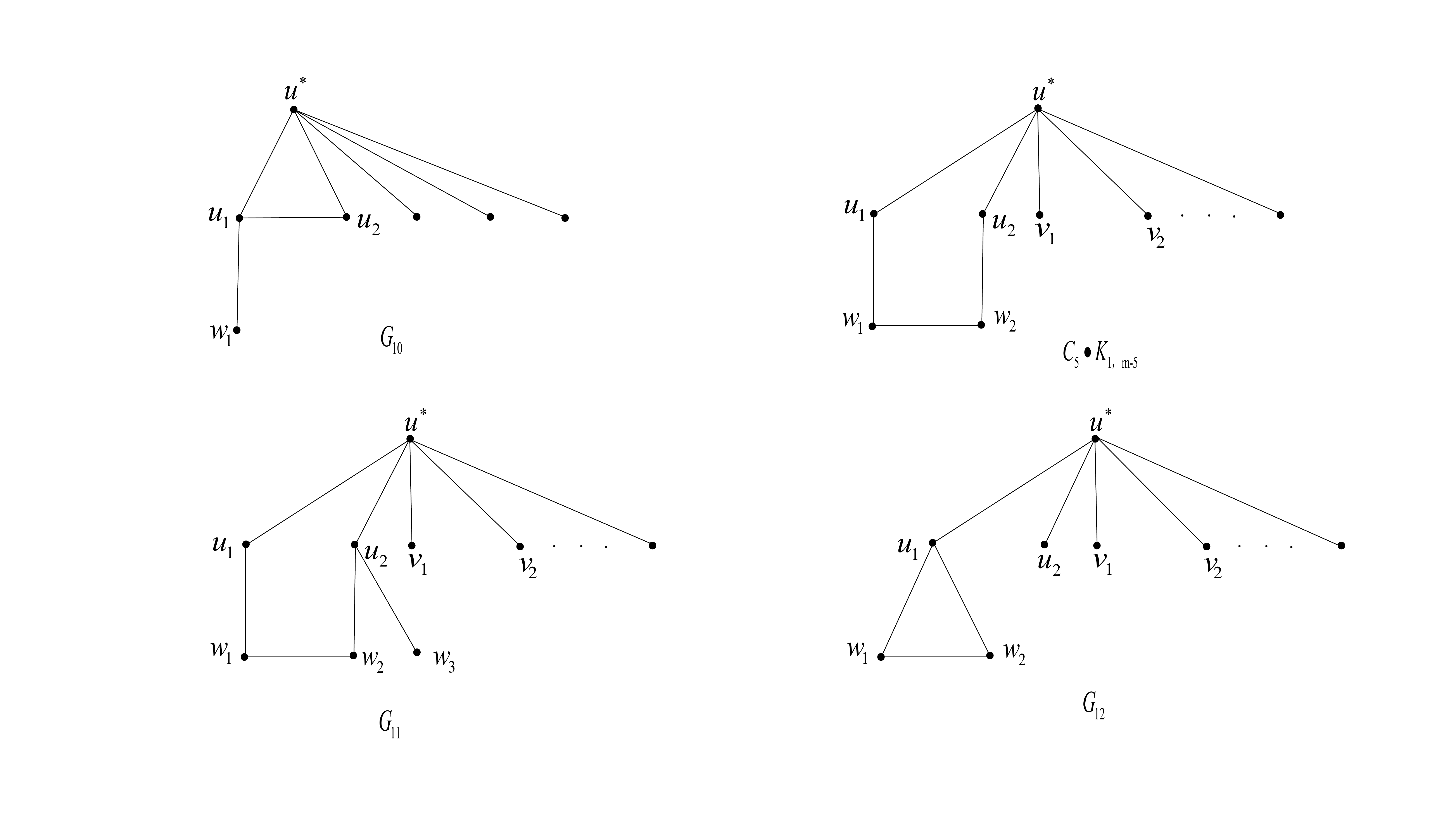}
\caption{Graphs $G_{10}-G_{12}$ and $C_{5}\bullet K_{1,m-5}$ of Subcases 2.1 and 2.2.}\label{3.}
\end{centering}
\end{figure}

\section{Proof of Theorem 1.6.}

Let $G^{\ast}$ be the extremal graph with maximum spectral radius in $\mathcal{G}(m,F)$ for a fixed $F$. Let $\rho^{\ast}=\rho(G^{\ast})$ and let $X^{\ast}$ be the Perron vector of $G^{\ast}$ with coordinate $x_{v}$ corresponding to the vertex $v\in V(G^{\ast})$. Recall that $W=V(G^{\ast})\backslash N[u^{\ast}]$. A vertex $u^{\ast}$ in $G^{\ast}$ is said to be an extremal vertex if $x_{u^{\ast}}=max\{x_{v}\mid v\in V(G^{\ast})\}$.
\noindent\begin{lemma}\label{le:ch-4.1.}{\rm(}$\cite{ZhLS}${\rm)} If $F$ is a $2$-connected graph and $u^{\ast}$ is an extremal vertex of  $G^{\ast}$, then the following statements hold.

{\rm(}$i${\rm)} $G^{\ast}$ is connected.

{\rm(}$ii${\rm)} There exists no cut vertex in $V(G^{\ast})\setminus\{u^{\ast}\}$ and hence $d(u)\geq2$ for any $u\in V(G^{\ast})\setminus N[u^{\ast}]$.

{\rm(}$iii${\rm)} If $F$ is $C_{4}$-free, then $N(u_{1})=N(u_{2})$ for any non-adjacent vertices of $u_{1}, u_{2}$ of degree two.
\end{lemma}
\noindent\begin{lemma}\label{le:ch-4.2.}{\rm(}$\cite{BrH}${\rm)} Let $G$ be a bipartite graph of size $m$. Then $\rho(G)\leq\sqrt{m}$, with equality if and only if $G$ is a disjoint union of a complete bipartite graph and isolated vertices.
\end{lemma}

\noindent\begin{lemma}\label{le:ch-4.3.}{\rm(}$\cite{MiLH}${\rm)} $\rho(S_{\frac{m+4}{2},2}^{-})>\frac{1+\sqrt{4m-5}}{2}$ for $m\geq6$.

\end{lemma}

\noindent\begin{lemma}\label{le:ch-4.4.}{\rm(}$\cite{MiLH}${\rm)} Let $X=\{x_{1}, x_{2}, \ldots, x_{n}\}^{T}$ be the Perron vector of a connected graph $G$ of size $m$ and let $x_{u^{\star}}=max\{x_{v}| v\in V(G)\}$. If $\rho(G)> \frac{1+\sqrt{4m-5}}{2}$, then we have the following results.

{\rm(}$i${\rm)}
\begin{equation}\label{eq:ch-10}
\sum_{v\in N(u^{\star})\setminus N_{0}(u^{\star})}(d_{N(u^{\star})}(v)-1)x_{v}>(e(W)+e(N(u^{\star}))-\frac{3}{2})x_{u^{\star}},
\end{equation}
and
\begin{equation}\label{eq:ch-11}
e(W)<e(N(u^{\star}))-|N(u^{\star})\setminus N_{0}(u^{\star})|+\frac{3}{2}
\end{equation}

{\rm(}$ii${\rm)} If there exists a vertex $v$ of $G$ such that $x_{v}< (1-\beta)x_{u^{\star}}$ where $0<\beta<1$, then
\begin{equation}\label{eq:ch-12}
e(W)<e(N(u^{\star}))-|N(u^{\star})\setminus N_{0}(u^{\star})|+\frac{3}{2}-\beta d_{N(u^{\star})}(v),  \mbox{for }  v\in N^{2}(u^{\star})\subseteq W,
\end{equation}
\begin{equation}\label{eq:ch-13}
e(W)<e(N(u^{\star}))-|N(u^{\star})\setminus N_{0}(u^{\star})|+\frac{3}{2}-\beta (d_{N(u^{\star})}(v)-1),  \mbox{for }  v\in N(u^{\star})\setminus N_{0}(u^{\star}).
\end{equation}
{\rm(}$iii${\rm)} If there exists a subset $S\subseteq N(u^{\star})\backslash N_{0}(u)$ such that $x_{v}< (1-\beta)x_{u^{\star}}$ for any $i\in V(S)$ and $0<\beta<1$, then
\begin{equation}\label{eq:ch-14}
e(W)<e(N(u^{\star}))-|N(u^{\star})\setminus N_{0}(u^{\star})|+\frac{3}{2}-\beta\sum_{v\in S}(d_{N(u^{\star})}(v)-1).
\end{equation}
\end{lemma}

\noindent\begin{lemma}\label{le:ch-4.5.} Let $G^{\ast}$ be a $C_{5}^{+}$-free graph with $u^{\ast}\in V(G)$ and $L$ be a component of $G^{\ast}[N(u^{\ast})]$. Then $L$ is one of the following statements.

{\rm(}$i${\rm)} a star $K_{1,r}$ for $r\geq0$, where $K_{1,0}$ is a singleton component.

{\rm(}$ii${\rm)} a double star $D_{a,b}$ for $a,b\geq1$.

{\rm(}$iii${\rm)} a copy of $S_{r+1}^{1}$ for $r\geq2$, where $S_{3}^{1}$ is a triangle for $r=2$.

{\rm(}$iv${\rm)} a graph with $C_{4}$ as its spanning subgraph, that is, $C_{4}$, $C_{3}^{+}$ or $K_{4}$.
\end{lemma}

\noindent {\bf Proof.} Since $G^{\ast}$ contains no $C_{5}^{+}$, then $G^{\ast}[N(u^{\ast})]$ contains no any path of length more than $3$ and any cycle of length more than $4$. If $G^{\ast}[N(u^{\ast})]$ contains $P_{1}$ as a subgraph, then $L\cong K_{1,0}$. If $G^{\ast}[N(u^{\ast})]$ contains $P_{2}$ as a subgraph, then $L\cong K_{1,1}$  or $L\cong K_{i}$ for each $i\in\{3,4\}$. If $G^{\ast}[N(u^{\ast})]$ contains $P_{3}$ as a subgraph, then $L\cong C_{3}^{+}, K_{1,r}$ or $S_{r+1}^{1}$ for $r\geq2$. If $G^{\ast}[N(u^{\ast})]$ contains $P_{4}$ as a subgraph, then $L\cong D_{a,b}$ for $a,b\geq1$, as desired. $\qedsymbol$

For each component $L$ of $G^{\ast}[N(u^{\ast})]$, let $W_{L}=\{w\mid w\in W\cap N_{u\in L}(u)\}$. Thus $W_{L_{i}}\cap W_{L_{j}}=\emptyset$ for any two distinct components $L_{i}$ and $L_{j}$ of $G^{\ast}[N(u^{\ast})]$, unless one of $L_{i}$ and $L_{j}$ is an isolated vertex and the other is a star $K_{1,r}$ for $r\geq0$ (that is, vertices in $W_{L_{i}}\cap W_{L_{j}}$ must be adjacent to the center vertex of the star $K_{1,r}$ for $r\geq0$).

Note that $\rho^{\ast}\geq \rho(S_{\frac{m+4}{2},2}^{-})> \frac{1+\sqrt{4m-5}}{2}>9$ for $m\geq74$. Thus ${\rho^{\ast}}^{2}-\rho^{\ast}>m-\frac{3}{2}$. Let $N_{+}(u^{\ast})=N(u^{\ast})\setminus N_{0}(u^{\ast})$. By \eqref{eq:ch-3}, we have
$$(m-\frac{3}{2})x_{u^{\ast}}<({\rho^{\ast}}^{2}-\rho^{\ast})x_{u^{\ast}}\leq |N(u^{\ast})|x_{u^{\ast}}+\sum_{v\in N_{+}(u^{\ast})}(d_{N(u^{\ast})}(v)-1)x_{v}+e(N(u^{\ast}), W)-\sum_{v\in N_{0}(u^{\ast})}x_{v}.$$
It follows that
$$\left(m-\frac{3}{2}-|N(u^{\ast})|-e(N(u^{\ast}), W)+\sum_{v\in N_{0}(u^{\ast})}\frac{x_{v}}{x_{u^{\ast}}}\right)x_{u^{\ast}}
<\sum_{v\in N_{+}(u^{\ast})}(d_{N(u^{\ast})}(v)-1)x_{v}.$$
Let $\zeta(L)=\sum_{v\in V(L)}(d_{L}(v)-1)x_{v}$. For each non-trivial connected component $L$ of $G^{\ast}[N(u^{\ast})]$,
we have
\begin{equation}\label{eq:ch-15}
\left(e(N(u^{\ast}))+e(W)+\sum_{v\in N_{0}(u^{\ast})}\frac{x_{v}}{x_{u^{\ast}}}-\frac{3}{2}\right)x_{u^{\ast}} <\sum_{L}\zeta(L).
\end{equation}
\noindent\begin{lemma}\label{le:ch-4.6.}
Let $G^{\ast}$ be the extremal graph which attains maximum spectral radius $\rho^{\ast}=\rho(G^{\ast})$ among all $C_{5}^{+}$-free graphs with even size $m\geq74$, and let $X=\{x_{1}, x_{2}, \ldots, x_{n}\}^{T}$ be the Perron vector of $G^{\ast}$ and $u^{\ast}$ be an extremal vertex. Let $L^{\ast}$ be a component of $G^{\ast}[N_{+}(u^{\ast})]$. If $\rho^{\ast}>\frac{1+\sqrt{4m-5}}{2}$, then

{\rm(}$i${\rm)} $G^{\ast}[N_{+}(u^{\ast})]$ does not contain $C_{4}$ as a spanning subgraph, that is, which does not contain one of $C_{4}$,$C_{3}^{+}$ and $K_{4}$ as a spanning subgraph.

{\rm(}$ii${\rm)} $e(W)=0$, furthermore, $L^{\ast}\ncong K_{3}$ for any component $L^{\ast}$ of $G^{\ast}[N_{+}(u^{\ast})]$.

{\rm(}$iii${\rm)} $G^{\ast}[N_{+}(u^{\ast})]$ has exactly one star component $K_{1,r}$ for some $r\geq3$ and $W=\emptyset$.
\end{lemma}

\noindent {\bf Proof.} {\rm(}$i${\rm)} Let $\mathcal{L}$ be the family of components of $G^{\ast}[N(u^{\ast})]$ each of which contains $C_{4}$ as a spanning subgraph and $\mathcal{L^{\prime}}$ be the family of other non-trivial components of $G^{\ast}[N(u^{\ast})]$ each of which contains no $C_{4}$ as a spanning subgraph. By Lemma \ref{le:ch-4.5.} (i)-(iii), for each $L\in\mathcal{L^{\prime}}$, we have
$$\zeta(L)=\sum_{v\in V(L)}(d_{L}(v)-1)x_{v}\leq(2e(L)-|V(L)|)x_{u^{\ast}}\leq e(L)x_{u^{\ast}}.$$ For any two distinct components $L_{i}, L_{j}\in\mathcal{L}$, since $G^{\ast}$ contains no $C_{5}^{+}$, we have $W_{L_{i}}\cap W_{L_{j}}=\emptyset$ and $e(W_{L_{i}}, W_{L_{j}})=0$. Hence, $e(W)\geq \sum_{L\in\mathcal{L}}e(W_{L}, W)$. By \eqref{eq:ch-15}, we have
\begin{equation}\label{eq:ch-16}
\left(\sum_{L\in \mathcal{L}}(e(L)+e(W_{L}, W))-\frac{3}{2}\right)x_{u^{\ast}}<\sum_{L\in\mathcal{L}}\zeta(L).
\end{equation}
Suppose that $\mathcal{L}\neq\emptyset$, we will show that $\zeta(L)\leq(e(L)+e(W_{L}, W))-\frac{3}{2})x_{u^{\ast}}$ holds for each $L\in \mathcal{L}$ and $\sum_{L\in\mathcal{L}}\zeta(L)\leq \left(\sum_{L\in \mathcal{L}}(e(L)+e(W_{L}, W)-\frac{3}{2})\right)x_{u^{\ast}}$ which contradicts \eqref{eq:ch-16}. Let $L^{\ast}\in \mathcal{L}$ with $V(L^{\ast})=\{u_{1},u_{2},u_{3},u_{4}\}$.

{\bf Case 1.}  $W_{L^{\ast}}=\emptyset$.

Assume that $x_{u_{1}}=max\{x_{u_{i}}: 1\leq i\leq4\}$. Hence, $\rho^{\ast}x_{u_{1}}=\sum_{u\in N(u_{1})}x_{u}\leq x_{u^{\ast}}+3x_{u_{1}}$, i.e., $x_{u_{1}}\leq\frac{x_{u^{\ast}}}{\rho^{\ast}-3}<\frac{x_{u^{\ast}}}{6}$ for $\rho^{\ast}>9$. Note that $4\leq e(L^{\ast})\leq6$. It follows that
$$ \zeta(L^{\ast})=\sum_{v\in V(L^{\ast})}(d_{L^{\ast}}(v)-1)x_{v}\leq(2e(L^{\ast})-4)x_{u_{1}}\leq \frac{1}{3}(e(L^{\ast})-2)x_{u^{\ast}}<(e(L^{\ast})-\frac{3}{2})x_{u^{\ast}},$$ as desired.

{\bf Case 2.} $W_{L^{\ast}}\neq\emptyset$.

Note that $d_{N(u^{\ast})}(w)=d_{L^{\ast}}(w)=1$ for $w\in W_{L^{\ast}}$. By Lemma \ref{le:ch-4.1.} (ii), we have $e(W_{L^{\ast}}, W)\geq1$. We consider the following three subcases.

{\bf Subcase 2.1.}  All vertices in $W_{L^{\ast}}$ have a unique common neighbor $u_{1}$, i.e., $N_{W}(u_{i})=\emptyset$ for each $i\in \{2,3,4\}$.

Assume that $x_{u_{2}}=max\{x_{u_{i}}: 2\leq i\leq 4\}$. Therefore, $$\rho^{\ast}x_{u_{2}}\leq x_{u^{\ast}}+x_{u_{1}}+x_{u_{3}}+x_{u_{4}}\leq2(x_{u^{\ast}}+x_{u_{2}}),$$ it follows that $x_{u_{2}}\leq \frac{2x_{u^{\ast}}}{\rho^{\ast}-2}<\frac{2x_{u^{\ast}}}{7}$. Note that $4\leq e(L^{\ast})\leq6$. Hence,
$$
\begin{aligned}
\zeta(L^{\ast})&=\sum_{v\in V(L^{\ast})}(d_{L^{\ast}}(v)-1)x_{v}\\
&\leq(d_{L^{\ast}}(u_{1})-1)x_{u_{1}}+(2e(L^{\ast})-d_{L^{\ast}}(u_{1})-3)x_{u_{2}}\\
&<(d_{L^{\ast}}(u_{1})-1)x_{u^{\ast}}+(\frac{4}{7}e(L^{\ast})-\frac{2}{7}d_{L^{\ast}}(u_{1})-\frac{6}{7})x_{u^{\ast}}\\
&\leq(\frac{4}{7}e(L^{\ast})+\frac{2}{7})x_{u^{\ast}}\\
&<(e(L^{\ast})-\frac{1}{2})x_{u^{\ast}}\\
&\leq (e(L^{\ast})+e(W_{L^{\ast}},W)-\frac{3}{2})x_{u^{\ast}},
\end{aligned}
$$
as desired.

{\bf Subcase 2.2.}  There exist exactly two vertices $w, w^{\prime}\in W_{L^{\ast}}$ with distinct neighbors in $V(L^{\ast})$.

In this case, we have $d(w)+d(w^{\prime})\geq5$ from Lemma \ref{le:ch-4.1.} (iii). By Lemma \ref{le:ch-4.1.} (ii), we have $e(W_{L^{\ast}}, W)=d(w)+d(w^{\prime})-e(\{w,w^{\prime} \},V(L^{\ast}))\geq3$. Since $L^{\ast}\in \mathcal{L}$, we have $e(L^{\ast})\leq6$ and $e(L^{\ast})\leq e(W_{L^{\ast}}, W)+3$. Let $N_{L^{\ast}}(w)=\{u_{1}\}, N_{L^{\ast}}(w^{\prime})=\{u_{2}\}, x_{u_{3}}=max\{x_{u_{3}}, x_{u_{4}}\}$, then $$\rho^{\ast}x_{u_{3}}\leq x_{u^{\ast}}+x_{u_{1}}+x_{u_{2}}+x_{u_{3}}\leq3x_{u^{\ast}}+x_{u_{3}},$$
and $$x_{u_{3}}\leq\frac{3}{\rho^{\ast}-1}<\frac{3}{8}x_{u^{\ast}}.$$ Since $d_{L^{\ast}}(u_{3}), d_{L^{\ast}}(u_{4})\geq2$, we obtain
$$
\begin{aligned}
\zeta(L^{\ast})&=\sum_{v\in V(L^{\ast})}(d_{L^{\ast}}(v)-1)x_{v}\leq \sum_{v\in V(L^{\ast})}(d_{L^{\ast}}(v)-1)x_{u^{\ast}}-2(d_{L^{\ast}}(u_{3})-1)(x_{u^{\ast}}-x_{u_{3}})\\
&\leq (2e(L^{\ast})-|L^{\ast}|)x_{u^{\ast}}-\frac{5}{4}x_{u^{\ast}}\leq(e(L^{\ast})+e(W_{L^{\ast}}, W)-\frac{9}{4})x_{u^{\ast}}\\
&<(e(L^{\ast})+e(W_{L^{\ast}}, W)-\frac{3}{2})x_{u^{\ast}},
\end{aligned}
$$ as desired.

{\bf Subcase 2.3.}  There exist $k$ ($k\geq3$) vertices, say $w_{1}, w_{2}, \ldots, w_{k}$, of $W_{L^{\ast}}$ such that they have mutual distinct neighbors in $V(L^{\ast})$.

In this case, if $w_{i}w_{j}\in E(G^{\ast}[W_{L^{\ast}}])$, then $N_{L^{\ast}}(w_{i})=N_{L^{\ast}}(w_{j})$. Hence, $\{w_{1}, w_{2}, \ldots, w_{k}\}$ is an independent set of $G^{\ast}$ from Lemma \ref{le:ch-4.1.} (iii). By Lemma \ref{le:ch-4.1.} (ii), we obtain that $d(w_{i})\geq2$ for $1\leq i\leq k$ and $d(w_{i})=2$ holds for at most one vertex $w_{i}$. Therefore, $\sum_{1\leq i\leq k}d(w_{i})\geq3k-1$ and $e(W_{L^{\ast}}, W)\geq2k-1$. Thus
$$ e(L^{\ast})\leq e(K_{4})=6\leq e(W_{L^{\ast}}, W)-2k+7\leq e(W_{L^{\ast}}, W)+1$$ and $\zeta(L^{\ast})=\sum_{v\in V(L^{\ast})}(d_{L^{\ast}}(v)-1)x_{v}\leq (2e(L^{\ast})-4)x_{u^{\ast}}\leq (e(L^{\ast})+e(W_{L^{\ast}}, W)-3)x_{u^{\ast}}<(e(L^{\ast})+e(W_{L^{\ast}}, W)-\frac{3}{2})x_{u^{\ast}},$ as desired. This completes the proof of {\rm(}$i${\rm)}.

{\rm(}$ii${\rm)} By Lemmas \ref{le:ch-4.5.} and \ref{le:ch-4.6.} (i), each component $L$ of $G[N_{+}(u^{\ast})]$ is either a tree or a unicyclic graph $S_{r+1}^{1}$ for some $r\geq2$. Let $\mathcal{L^{\prime}}$ be the family of the components of $G[N_{+}(u^{\ast})]$. Assume that there are $c$ non-trivial tree components in $G^{\ast}[N_{+}(u^{\ast})]$, then
$$\sum_{L\in\mathcal{L^{\prime}}}\zeta(L)=\sum_{L\in\mathcal{L^{\prime}}}\sum_{v\in V(L)}(d_{L}(v)-1)x_{v}\leq \sum_{L\in\mathcal{L^{\prime}}}(2e(L)-|V(L)|)x_{u^{\ast}}=\left(e(N(u^{\ast}))-c \right)x_{u^{\ast}},$$
where $L\in \mathcal{L^{\prime}}$ takes over all non-trivial components of $G^{\ast}[N_{+}(u^{\ast})]$. Combining with \eqref{eq:ch-15}, we have
\begin{equation}\label{eq:ch-17}
e(W)< \frac{3}{2}-c-\sum_{v\in N_{0}(u^{\ast})}\frac{x_{v}}{x_{u^{\ast}}}.
\end{equation}
Hence, $e(W)\leq1$ and $c\leq1$. In addition, $e(W)=1$ holds if and only if $c=0$ and $\sum_{v\in N_{0}(u^{\ast})}\frac{x_{v}}{x_{u^{\ast}}}<\frac{1}{2}$. Then each component $L$ of $G^{\ast}[N_{+}(u^{\ast})]$ contains $C_{3}$ as a subgraph. Without loss of generality, let $w_{1}w_{2}$ be the unique edge in $E(W)$. If $\{w_{1}, w_{2}\}\in N_{W}(L)$, then there exists a cut vertex or $C_{5}^{+}$ as a subgraph. If $w_{1}\in N_{W}(L), w_{2}\in N_{W}(N_{0}(u^{\ast}))$, then there exists $C_{5}^{+}$ as a subgraph. Thus $\{w_{1}, w_{2}\}\in N_{W}(N_{0}(u^{\ast}))$. By Lemma \ref{le:ch-4.1.} (ii), $d_{N_{0}(u^{\ast})}(w_{i})\geq1$ for each $i\in\{1,2\}$.

We claim that $|N_{N_{0}(u^{\ast})}(w_{1})\cap N_{N_{0}(u^{\ast})}(w_{2})|\leq2$, otherwise, there exists $C_{5}^{+}$ as a subgraph. Let $x_{w_{1}}=max\{x_{w_{1}}, x_{w_{2}}\}$. $$\rho^{\ast}x_{w_{1}}=x_{w_{2}}+\sum_{v\in N_{N_{0}(u^{\ast})}(w_{1})}x_{v}\leq x_{w_{1}}+\sum_{v\in N_{0}(u^{\ast})}x_{v}<x_{w_{1}}+\frac{1}{2}x_{u^{\ast}},$$ i.e., $$x_{w_{1}}<\frac{1}{2(\rho^{\ast}-1)}x_{u^{\ast}}<\frac{1}{16}x_{u^{\ast}}.$$ By \eqref{eq:ch-12}, we have $$e(W)<\frac{3}{2}-\frac{15}{16}d_{N_{0}(u^{\ast})}(w_{1})\leq\frac{9}{16},$$ a contradiction. Thus $e(W)=0$. By \eqref{eq:ch-17}, we have $$\frac{3}{2}-c-\sum_{v\in N_{0}(u^{\ast})}\frac{x_{v}}{x_{u^{\ast}}}>0.$$ Furthermore, we have either $c=0$ and
\begin{equation}\label{eq:ch-18}
\sum_{v\in N_{0}(u^{\ast})}\frac{x_{v}}{x_{u^{\ast}}}<\frac{3}{2}
\end{equation}
or $c=1$ and
\begin{equation}\label{eq:ch-19}
 \sum_{v\in N_{0}(u^{\ast})}\frac{x_{v}}{x_{u^{\ast}}}<\frac{1}{2}.
\end{equation}
If $c=0$ and $\sum_{v\in N_{0}(u^{\ast})}\frac{x_{v}}{x_{u^{\ast}}}<\frac{3}{2}$, then $G^{\ast}[N_{+}(u^{\ast})]$ contains a component $L^{\ast}\cong S_{r+1}^{1}$ for some $r\geq2$.

Suppose that $L^{\ast}\cong K_{3}$ with $V(L^{\ast})=\{u_{1}, u_{2}, u_{3}\}$. If $W_{L^{\ast}}=\emptyset$, then $$x_{u_{1}}=x_{u_{2}}=x_{u_{3}}=\frac{x_{u^{\ast}}}{\rho^{\ast}-2}<\frac{x_{u^{\ast}}}{7}.$$ Hence, $$\zeta(L^{\ast})=\sum_{1\leq i\leq3}(d_{L^{\ast}}(u_{i})-1)x_{u_{i}}=3x_{u_{1}}<\frac{3}{7}x_{u^{\ast}}=\frac{3}{7}(e(L^{\ast})-2)x_{u^{\ast}}.$$
Since $e(W)=0$ and $\zeta(L)\leq e(L)x_{u^{\ast}}$ for each non-trivial component $L\in \mathcal{L^{\prime}}\backslash L^{\ast}$ of $G^{\ast}[N_{+}(u^{\ast})]$. Combining with \eqref{eq:ch-15}, we have
$$\left(e(N(u^{\ast}))+\sum_{v\in N_{0}(u^{\ast})}\frac{x_{v}}{x_{u^{\ast}}}-\frac{3}{2}\right)x_{u^{\ast}} <\sum_{L\in\mathcal{L^{\prime}}}\zeta(L)<\left(e(N(u^{\ast}))-\frac{4}{7}e(L^{\ast})-\frac{6}{7}\right)x_{u^{\ast}},$$ it follows that $\sum_{v\in N_{0}(u^{\ast})}\frac{x_{v}}{x_{u^{\ast}}}<\frac{-15}{14}$, a contradiction.
Thus $W_{L^{\ast}}\neq\emptyset$. Note that $2\leq d(w)\leq V(L^{\ast})=3$ for each vertex $w\in W_{L^{\ast}}$ and $N_{W}(L^{\ast})\cap N_{W}(N_{0}(u^{\ast}))=\emptyset$. If there is a vertex $w\in W_{L^{\ast}}$ such that $d(w)=3$, then $W_{L^{\ast}}=\{w\}$. If $d(w)=2$ for each vertex $w\in W_{L^{\ast}}$, then Lemma \ref{le:ch-4.1.} (iii) implies that all vertices in $W_{L^{\ast}}$ share the same neighborhoods. Without loss of generality, assume that $N(w)=\{u_{1}, u_{2}\}$ for each vertex $w\in W_{L^{\ast}}$. Let $G_{13}=G^{\ast}-\{wu_{1}|w\in N_{W_{L^{\ast}}}(u_{1}) \}+\{wu^{\ast}|w\in N_{W_{L^{\ast}}}(u_{1})\}$. In both cases, we have $G_{13}\in \mathcal{G}(m, C_{5}^{+})$ and $\rho(G_{13})>\rho^{\ast}$ from Lemma \ref{le:ch-2.1.}, a contradiction. Thus $G^{\ast}[N_{+}(u^{\ast})]$ contains no a component $L^{\ast}\cong K_{3}$. This completes the proof of {\rm(}$ii${\rm)}.

{\rm(}$iii${\rm)} Suppose that $L^{\ast}\cong S_{r+1}^{1}$ for some $r\geq3$, then we will prove that $L^{\ast}$ is the unique non-trivial component of $G^{\ast}[N_{+}(u^{\ast})]$. Note that there are $r-2$ vertices in $V(L^{\ast})$ with degree two in $G^{\ast}$. By Lemma \ref{le:ch-4.1.} (iii), there does not exist a vertex of degree two out of other components. Then $L^{\ast}$ is the unique component which contains $K_{3}$ as a subgraph. In this case, we suppose that $G^{\ast}[N_{+}(u^{\ast})]$ contains another non-trivial tree component $L$. By Lemma \ref{le:ch-4.1.} (ii), we obtain that $d(w)\geq2$ for each $w\in W$. Combining with $e(W)=0$, we obtain that $W_{L^{\ast}}=\emptyset$ and $L$ is a tree. In addition, $W_{L}\neq\emptyset$ and $d(w)\geq3$ for each vertex $w\in W_{L}\cup V(L)$. Since $e(W)=0$ and $W_{L}\cap W_{L^{\ast}}=\emptyset$. Then $N(w)\subseteq V(L)$ for each vertex $w\in W_{L}$. Let $V(L^{\ast})=\{u_{0}, u_{1}, \ldots, u_{r}\}$ with $d_{L^{\ast}}(u_{0})=r$ and $d_{L^{\ast}}(u_{1})=d_{L^{\ast}}(u_{2})=2$. Thus $x_{u_{1}}=x_{u_{2}}$ and $x_{u_{3}}=x_{u_{4}}=\ldots= x_{u_{r}}$.
Note that
$$\rho^{\ast}x_{u_{1}}= x_{u_{0}}+x_{u_{2}}+x_{u^{\ast}}\leq x_{u_{1}}+2x_{u^{\ast}}.$$

It follows that $$x_{u_{1}}\leq \frac{2x_{u^{\ast}}}{\rho^{\ast}-1}<\frac{x_{u^{\ast}}}{4}$$ for $\rho^{\ast}>9$.
By \eqref{eq:ch-14}, $$e(W)<\frac{3}{2}-1-\frac{3}{4}\sum_{v\in\{u_{1},u_{2}\}}(d_{N(u^{\ast})}(v)-1)=-1,$$ a contradiction. Thus there is a non-trivial unique component of $G^{\ast}[N_{+}(u^{\ast})]$. Since
$$\zeta(L^{\ast})=(r-1)x_{u_{0}}+x_{u_{1}}+x_{u_{2}}\leq (r-1)x_{u^{\ast}}+2x_{u_{1}}<(r-\frac{1}{2})x_{u^{\ast}}=(e(L^{\ast})-\frac{3}{2})x_{u^{\ast}}.$$
Combining with $e(W)=0$ and \eqref{eq:ch-15}, we have
$$(e(N(u^{\ast}))+\sum_{v\in N_{0}(u^{\ast})}\frac{x_{v}}{x_{u^{\ast}}}-\frac{3}{2})x_{u^{\ast}}     <\zeta(L^{\ast})<(e(N(u^{\ast}))-\frac{3}{2})x_{u^{\ast}},$$ it follows that $\sum_{v\in N_{0}(u^{\ast})}\frac{x_{v}}{x_{u^{\ast}}}<0$, a contradiction. Hence,  $G^{\ast}[N_{+}(u^{\ast})]$ contains no unicyclic graph and contains $c$ non-trivial tree components. If $c=0$, then $G^{\ast}$ is bipartite. By Lemma \ref{le:ch-4.2.}, we have $\rho^{\ast}\leq\sqrt{m}< \frac{1+\sqrt{4m-3}}{2}$ for $m\geq74$, a contradiction. Thus $c=1$ and \eqref{eq:ch-19} holds, i.e., $G^{\ast}[N_{+}(u^{\ast})]\cong L$, where $L$ is a non-trivial tree. By Lemma \ref{le:ch-4.5.}, $diam(L)\leq3$.

If $diam(L)\leq3$, then $L$ is a double star. Since $G^{\ast}$ is $C_{5}^{+}$-free, we have $d_{N(u^{\ast})}(w)=1$ for each vertex $w\in W_{L}$. Combining with $e(W)=0$ and Lemma \ref{le:ch-4.1.} (ii), we have $W_{L}=\emptyset$, then $G^{\ast}$ contains two non-adjacent vertices of degree two with distinct neighborhoods, which contradicts the Lemma \ref{le:ch-4.1.} (iii).
Thus $diam(L)\leq2$, then $L\cong K_{1,r}$ for some $r\geq1$.

Let $V(L)=\{u_{0}, u_{1}, \ldots, u_{r}\}$ with center vertex $u_{0}$ and $d_{L}(u_{0})=r\geq1$. By Lemma \ref{le:ch-4.1.} (ii), we have $d_{N(u^{\ast})}(w)\geq2$ for any vertex $w\in W$. For $r=1$, we have $9x_{u^{\ast}}<\rho^{\ast}x_{u^{\ast}}=x_{u_{0}}+x_{u_{1}}+\sum_{v\in N_{0}(u^{\ast})}x_{v}<\frac{5}{2}x_{u^{\ast}}$. For $r=2$, we have $9x_{u^{\ast}} <\rho^{\ast}x_{u^{\ast}} =x_{u_{0}} +x_{u_{1}}+x_{u_{2}}+\sum_{v\in N_{0}(u^{\ast})}x_{v}<\frac{7}{2}x_{u^{\ast}}$, a contradiction. For $r\geq3$, we discuss the following three cases.

{\bf Case 1.}  $d_{L}(w)=1$.

 In this case, we obtain that $w$ is only adjacent to the center vertex $u_{0}$. By Lemma \ref{le:ch-4.2.} (iii), we have $d_{N_{0}(u^{\ast})}(w)\geq2$. $\rho^{\ast}x_{w}=x_{u_{0}}+\sum_{v\in N_{N_{0}(u^{\ast})}(W)}x_{v}\leq x_{u^{\ast}}+\frac{1}{2}x_{u^{\ast}}=\frac{3}{2}x_{u^{\ast}}$, i.e., $x_{w}\leq \frac{3}{2\rho^{\ast}}x_{u^{\ast}}<\frac{1}{6}x_{u^{\ast}}$. By \eqref{eq:ch-12}, we have $$e(W)<e(N(u^{\ast}))-|N_{+}(u^{\ast})|+\frac{3}{2}-\frac{5}{6}d_{N(u^{\ast})}(w)=\frac{1}{2}-\frac{5}{2}<0,$$ a contradiction.

{\bf Case 2.}  $d_{L}(w)=2$.

In this case, we have $N_{N_{0}(u^{\ast})}(w)=\emptyset$, otherwise, there is $C_{5}^{+}$ as a subgraph. By Lemma \ref{le:ch-4.2.} (iii), there exist two non-adjacent vertices of degree two with distinct neighborhoods, a contradiction.

{\bf Case 3.}  $d_{L}(w)\geq3$.

In this case, we have $N_{N_{0}(u^{\ast})}(w)=\emptyset$, otherwise, there is $C_{5}^{+}$ as a subgraph. Let $\{u_{0}, u_{1}, u_{2}\}\in N_{L}(w)$. Thus $G^{\ast}[u^{\ast},u_{1},w,u_{2},u_{0},u_{3}]$ contains $C_{5}^{+}$ as a subgraph, a contradiction. Let $\{u_{1}, u_{2}, u_{3}\}\in N_{L}(w)$ and $N_{N_{0}(u^{\ast})}(w)=\emptyset$. Thus $G^{\ast}[u^{\ast},u_{1},w,u_{2},u_{0},u_{3}]$ contains $C_{5}^{+}$ as a subgraph, a contradiction.

 By Case 1-3, we have $W_{L}=\emptyset$. By Lemma \ref{le:ch-4.6.} (ii), we have $e(W)=0$. Suppose that $W\neq \emptyset$, by Lemma \ref{le:ch-2.1.} (i), we obtain that $G^{\ast}$ is a connected graph. Thus $d(w)=d_{N_{0}(u^{\ast})}(w)$ for any vertex $w\in W$, furthermore, $N_{0}(u^{\ast})\neq \emptyset$. Combining with \eqref{eq:ch-19}, we have

 $$\rho^{\ast}x_{w}=\sum_{v\in N(w)} x_{v}\leq \sum_{v\in N_{0}(u^{\ast})}x_{v}<\frac{1}{2}x_{u^{\ast}},$$

it follows that $x_{w}< \frac{x_{u^{\ast}}}{2\rho^{\ast}}<\frac{x_{u^{\ast}}}{18}$.
By \eqref{eq:ch-12}, we have
$$e(W)<e(N(u^{\ast}))-|N_{+}(u^{\ast})|+\frac{3}{2}- \frac{17}{18}d_{N_{0}(u^{\ast})}(w)\leq\frac{1}{2}-\frac{17}{6}<0,$$
a contradiction. Thus $W=\emptyset$. $\qedsymbol$

\noindent\begin{lemma}\label{le:ch-4.7.} $G^{\ast}\cong S_{\frac{m+4}{2},2}^{-}$.
\end{lemma}

\noindent {\bf Proof.} By Lemmas \ref{le:ch-4.6.}, we have $e(W)=0, W=\emptyset$ and $G^{\ast}[N_{+}(u^{\ast})]\cong K_{1,r}$ for $r\geq3$. Thus $G^{\ast}\cong G_{14}$ (see Figure. 4). Let $|N_{0}(u^{\ast})|=t$. Since $m$ is even, we obtain that $t$ is odd and $t\geq1$. By Lemma \ref{le:ch-2.4.}, we obtain that $\rho^{\ast}$ is the largest root of the equation $f(x,t)=0$ where $$f(x,t)=x^{4}-mx^{2}-(m-t-1)x+\frac{t(m-t-1)}{2}$$for $m=t+1+2r\geq74$. Since $$f(x,t)-f(x,1)=(t-1)x+\frac{m(t-1)-t^{2}-t+2}{2}>0$$ for $x>0$ and $t\geq3$, this implies that $t=1$ for the extremal graph $G^{\ast}$. By Lemma \ref{le:ch-4.3.}, we have $\rho(S_{\frac{m+4}{2},2}^{-})> \frac{1+\sqrt{4m-5}}{2}$ for $m\geq74$ and $G^{\ast}\cong S_{\frac{m+4}{2},2}^{-}$, as desired.

This completes the proof of Theorem 1.6. $\blacksquare$
\begin{figure}[H]
\begin{centering}
\includegraphics[scale=0.25]{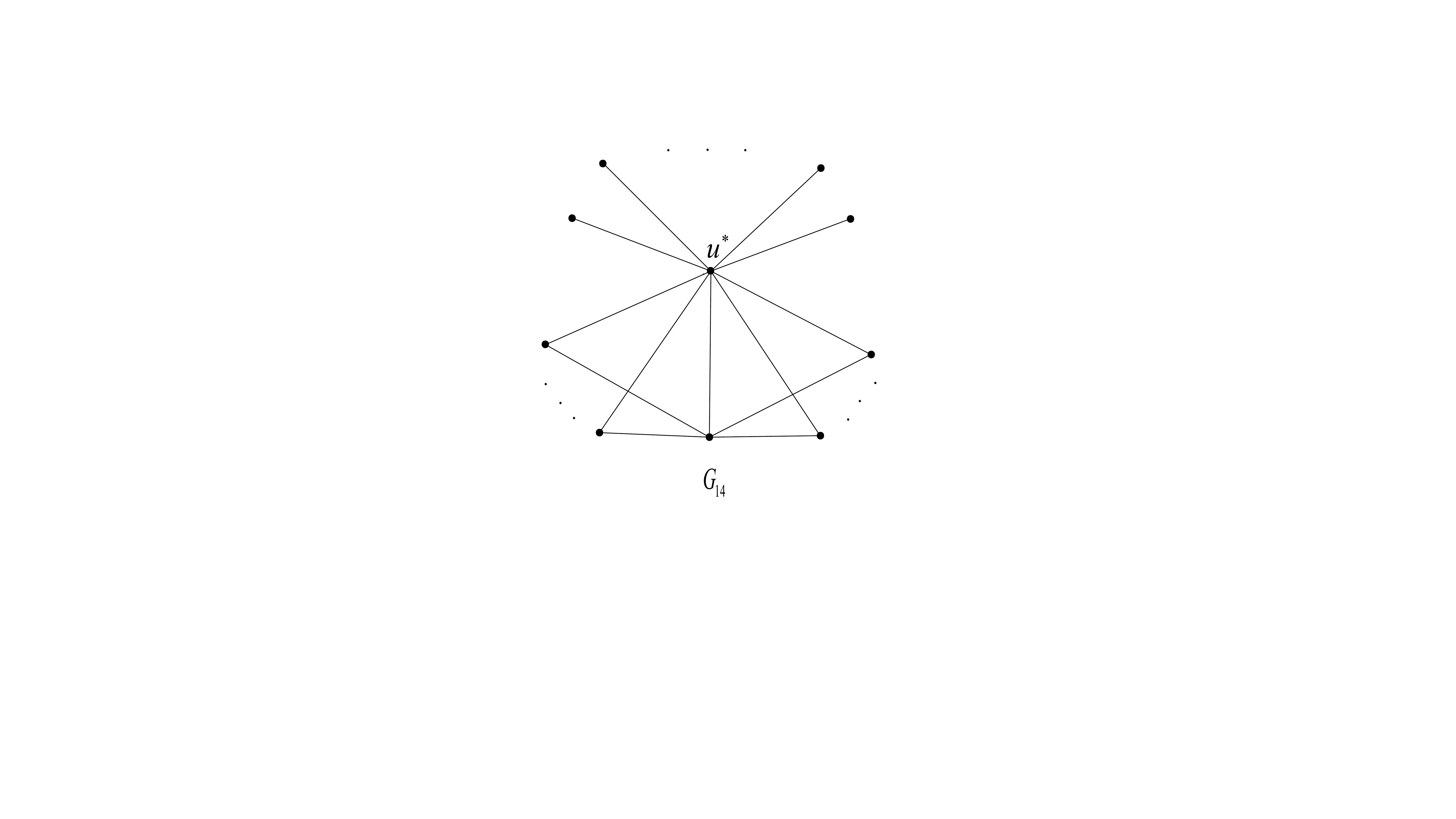}
\caption{Graph $G_{14}$ of Lemma 4.7.}\label{fi:ch-4}
\end{centering}
\end{figure}

\end{document}